\documentclass[12pt,reqno]{amsart} %
\usepackage{amsmath,amscd}
\usepackage{amsthm}
\usepackage{amssymb}
\usepackage{amsfonts}
\usepackage{latexsym}
\usepackage{url}
\usepackage{latexsym}
\usepackage{graphicx,color,subfig}
\usepackage{enumerate,fancybox}

\DeclareMathOperator\dist{dist}
\DeclareMathOperator\cpct{cap}

\DeclareMathOperator{\supp}{\mathrm{supp}}

\newtheorem{theorem}{Theorem}[section]

\newtheorem{definition}{Definition}[section]

\newtheorem{corollary}{Corollary}[section]

\numberwithin{equation}{section}

\def\proof#1. {\par
                      \ifdim\lastskip<15pt
                      \removelastskip\penalty-200
                      \vskip15pt plus3pt minus3pt
                      \fi
                       {\def\a{#1}
                       \ifx\a\empty
                       {\noindent\bf Proof.}
                       \else
                       {\noindent\bf Proof of #1.}
                       \fi}\enspace}
\def\restr#1{\,\vrule\,\lower1.75ex\hbox{$#1$}}

\input colordvi

\begin{document}
\title[]
{On the zeros of asymptotically extremal polynomial sequences in the plane}

\date{\today}

\thanks{}

\thanks{{\it Acknowledgements.}
The first author was partially supported  by the
U.S. National Science Foundation grant DMS-1109266.
The second author was supported by the University of Cyprus grant 3/311-21027.}

\author[Edward B.\ Saff]{E.B.\ Saff}
\address{Center for Constructive Approximation,
Department of Mathematics, Vanderbilt University,
1326 Stevenson Center, 37240 Nashville, TN, USA}
\email{edward.b.saff@Vanderbilt.Edu}
\urladdr{https://my.vanderbilt.edu/edsaff/}

\author[Nikos Stylianopoulos]{N. Stylianopoulos}
\address{Department of Mathematics and Statistics,
         University of Cyprus, P.O. Box 20537, 1678 Nicosia, Cyprus}
\email{nikos@ucy.ac.cy}
\urladdr{http://ucy.ac.cy/\textasciitilde nikos}


\begin{abstract}
Let $E$ be a compact set of positive logarithmic capacity  in the complex plane and let
$\{P_n(z)\}_{1}^{\infty}$ be a sequence of asymptotically extremal monic polynomials for $E$
 in the sense that
\begin{equation*}
\limsup_{n\to\infty}\|P_n\|_E^{1/n}\le\mathrm{cap}(E).
\end{equation*}
The purpose of this note is to provide sufficient geometric conditions on $E$ under which the (full) sequence of normalized counting measures
of the zeros of $\{P_n\}$ converges in the weak-star topology to the equilibrium measure on $E$, as $n\to\infty.$ Utilizing an argument of Gardiner and Pommerenke
dealing with the balayage of measures, we show that this is true, for example, if the interior of the polynomial convex hull of $E$ has a single component and the boundary
of this component has an ``inward corner" (more generally, a ``non-convex singularity"). This simple fact
 has thus far not been sufficiently
emphasized in the literature.
As applications we mention improvements of some known results on the distribution of zeros of
some special polynomial sequences.
\end{abstract}

\maketitle
\allowdisplaybreaks
\textit{Dedication:} To Herbert Stahl, an exceptional mathematician, a delightful personality, and a dear friend.
\section{Introduction}
Let $E$ be a compact set of positive logarithmic capacity (\textrm{cap}$(E)>0$) contained in the complex plane $\mathbb{C}$.
We denote  by $\Omega$ the unbounded component of $\overline{\mathbb{C}}\setminus E$ and
by $\mu_E$ the \textit{equilibrium measure} (energy minimizing Borel probability measure on $E$)
for the logarithmic potential on $E$; see e.g. \cite[Ch.~3]{Ra} and
\cite[Sect. I.1]{ST}. As is well-known, the support $\mathrm{supp}(\mu_E)$ lies on the boundary $\partial\Omega$ of $\Omega$.

For any polynomial $p_n(z)$,  of degree $n$, we denote by $\nu_{p_n}$ the \textit{normalized counting measure} for the zeros of $p_n(z)$; that is,
\begin{equation}
\nu_{p_n}:=\frac{1}{n}\sum_{p_n(z)=0}\delta_z,
\end{equation}
where $\delta_z$ is the unit point mass (Dirac delta) at the point $z$.

Let $\mathcal{N}$ denote an increasing sequence of positive integers. Then,
following \cite[p.~169]{ST} we say that a sequence of monic polynomials $\{P_n(z)\}_{n\in\mathcal{N}}$, of respective degrees $n$, is
\textit{asymptotically extremal on}  $E$ if
\begin{equation}
\limsup_{n\to\infty,\,n\in\mathcal{N}}\|P_n\|_E^{1/n}\le\mathrm{cap}(E),
\end{equation}
where $\|\cdot\|_E$ denotes the uniform norm on $E$.
(We remark that this inequality implies equality for the limit, since $\|P_n\|_E\ge \mathrm{cap}(E)^n$).
Such sequences arise, for example, in the study of  polynomials orthogonal with respect
to a measure $\mu$ belonging to the class \textbf{Reg}, see Definition 3.1.2 in \cite{StTobo}.

Concerning the asymptotic behavior of the zeros of an asymptotically extremal sequence of polynomials, we recall
the following result, see e.g. \cite[Thm 2.3]{MhSa91} and \cite[Thm III.4.7]{ST}.
\begin{theorem}\label{th:ST4.7}
Let $\{P_n\}_{n\in\mathcal{N}}$, denote an asymptotically extremal sequence of monic polynomials on $E$.
If $\mu$ is any weak$^*$ limit measure of the sequence $\{\nu_{P_n}\}_{n\in\mathcal{N}}$, then $\mu$ is a Borel probability
measure supported on $\overline{\mathbb{C}}\setminus\Omega$ and
$\mu^b=\mu_E$, where $\mu^b$ is the balayage of $\mu$ out of $\overline{\mathbb{C}}\setminus\Omega$ onto $\partial\Omega$.
Similarly, the sequence of balayaged counting measures converges to $\mu_E$:
\begin{equation}
\nu^b_{P_n}\,{\stackrel{*}{\longrightarrow}}\, \mu_E,\quad n\to\infty,\quad n\in\mathcal{N}.
\end{equation}
\end{theorem}
By the weak$^*$ convergence of a sequence of measures $\lambda_n$ to a measure $\lambda$ we mean that, for any continuous
$f$ with compact support in $\mathbb{C},$ there holds
\begin{equation*}
\int f d\lambda_n \to \int fd\lambda, \quad\textup{as }n\to\infty.
\end{equation*}
For properties of balayage,  see \cite[Sect. II.4]{ST}.

The goal of the present paper is to describe simple geometric conditions under which
the normalized counting measures $\nu_{P_n}$ of an asymptotically extremal sequence $\{P_n\}$  on $E$,
\textit{themselves} converge weak$^*$ to the equilibrium measure. For example, this is the case
whenever $E$ is a non-convex polygonal region, a simple fact that has thus far not been sufficiently emphasized in the
literature. Here we introduce  more general sufficient conditions based on arguments of Gardiner and Pommerenke \cite{GaPo02}
dealing with the balayage of measures.

The outline of the paper is as follows: In Section~2 we describe a geometric condition, which we call
the \textit{non-convex singularity} (NCS) condition and state the main result regarding the counting measures
$\nu_{P_n}$ of the zeros of  polynomials that form an asymptotically extremal sequence. Its proof is given in Section~4.

In Section~3, we apply the main result to obtain improvements in several previous results on the behavior of
the zeros of orthogonal polynomials, whereby the NCS condition yields convergence conclusions for the full
sequence $\mathcal{N}$ rather than for some subsequence.

\section{A geometric property}

\begin{definition}\label{Def1}
Let $G$ be a bounded simply connected domain in the complex plane. A point $z_0$ on the boundary
of $G$ is said to be a \textbf{ non-convex type singularity} (NCS) if it satisfies the following two conditions:
\begin{itemize}
\item[(i)]
There exists a closed disk $\overline{D}$ with $z_0$ on its circumference, such that $\overline{D}$ is contained in $G$ except for the
point $z_0$.
\item[(ii)]
There exists a line segment $L$ connecting a point $\zeta_0$ in the interior of $\overline{D}$ to $z_0$ such that
\begin{equation}\label{eq:Def1(ii)}
\lim_{z\to z_0 \atop{z\in L}}\frac{g_{G}(z,\zeta_0)}{|z-z_0|}=+\infty,
\end{equation}
where $g_{G}(z,\zeta_0)$ denotes the Green function of $G$ with pole at $\zeta_0\in G$.
\end{itemize}
\end{definition}
Recall that $g_G(\cdot,\zeta)$ is a positive harmonic function in $G\setminus\{\zeta\}$.
Also note that the assumption that $G$ is bounded and simply connected
implies that $G$ is regular with respect to the Dirichlet problem in $G$. This means that
$\lim_{z\to t \atop{z\in G}}g_G(z,\zeta)=0$, for any $t\in\partial G$; see, e.g.,
\cite[pp. 92 and 111]{Ra}.\\

\noindent\textbf{Remark.}
With respect to condition (ii), we note that the existence of a straight line $L$ and a point $\zeta_0$ for which
(\ref{eq:Def1(ii)}) holds, implies that the same is true for any other straight line connecting a point in the open disk
${D}$ with $z_0$. This can be easily deduced from Harnack's Lemma (see e.g. \cite [Lemma 4.9, p.17]{ST},\,\cite[p. 14]{ArGa}) in conjunction
with the symmetry property of Green functions, which imply that for a compact set $K \subset G$ containing $\zeta_0$ and $\dist (z,K)>\delta>0,\, z \in G$, there
 is constant $C=C(K,\delta, \zeta_0)>0$ such that the
inequality
\begin{equation}\label{eq:Har}
g_G(z,\zeta)\ge C g_G(z,\zeta_0)
\end{equation}
holds for all $\zeta \in K$.\\

As we shall easily show, a point $z_0$ satisfying the following sector condition is an NCS point.
\begin{definition}\label{Def2}
Let $G$ be a bounded simply connected domain. A point $z_0$ on the boundary of $G$ is said to be an \textbf{inward-corner}
(IC) \textbf{point} if there exists a circular sector of the form
$S:=\{z:0< |z-z_0|< r,\,\alpha\pi<\textup{arg}(z-z_0)<\beta\pi\}$
with $\beta-\alpha > 1$  whose closure is
contained in $G$ except for $z_0$.
\end{definition}
To see that an inner-corner point satisfies Definition \ref{Def1}, let $g_S(z,\zeta)$ denote the Green function of the sector $S$.
Then $g_S(z,\zeta)=-\log|\varphi_\zeta(z)|$, where $\varphi_\zeta$ is a conformal mapping of $S$ onto the unit disc $\mathbb{D}:=\{z:|z|<1\}$, satisfying $\varphi_\zeta(\zeta)=0$.
From the theory of conformal mapping it is known \cite{Lehman} that the following expansion
is valid for any $z\in\overline{S}$ near $z_0$:
\begin{equation*}
\varphi_\zeta(z)=\varphi_\zeta(z_0)+a_1(z-z_0)^{1/(\beta-\alpha)}+O(|z-z_0|^{2/(\beta-\alpha)}),
\end{equation*}
with $a_1\ne 0$. Since $|\varphi_\zeta(z_0)|=1$, the above implies that the limit in (2.1) holds with $g_S(z,\zeta_0)$ in the place
of $g_G(z,\zeta_0)$.  The desired limit then follows from the comparison principle for Green functions:
\begin{equation*}
g_S(z,\zeta)\le g_G(z,\zeta),\quad z,\zeta\in S;
\end{equation*}
see, e.g., \cite[p. 108]{Ra}.\\

\noindent\textbf{Remark.} It is interesting to note that if the boundary $\partial G$ is a piecewise analytic Jordan curve, then at any IC point $z_0$
of $\partial G$ the density of the equilibrium measure is zero. This can be easily deduced from the relation
$d\mu_{\partial G}(z)=|\Phi^\prime(z)|ds$ connecting the equilibrium measure to the arclength measure $ds$
on $\partial G$, where $\Phi$ is a conformal mapping of $\Omega$ onto $\{w:|w|>1\}$, taking $\infty$ to $\infty$.
Then, if $\lambda\pi$ ($1<\lambda<2$) is the interior opening angle at $z_0$, $\Phi(z)$ admits has near $z_0$
an expansion of the form
\begin{equation}\label{eq:Phi-exp}
\Phi(z)=\Phi(z_0)+b_1(z-z_0)^{1/(2-\lambda)}+o(|z-z_0|^{1/(2-\lambda)}),
\end{equation}
with $b_1\ne 0$, which leads to $\Phi^\prime(z_0)=0$.\\

We can now state our main result.
\begin{theorem}\label{th:main}
Let $E\subset\mathbb{C}$ be a compact set of positive capacity,  $\Omega$ the unbounded component
of $\overline{\mathbb{C}}\setminus E$, and $\mathcal{E}:=\overline{\mathbb{C}}\setminus{\Omega}$ denote the polynomial
convex hull of $E$. Assume there is closed set $E_0\subset\mathcal{E}$ with the following three properties:
\begin{itemize}
\item[(i)] {\rm{cap}}$(E_0)>0$;
\item[(ii)]
either $E_0=\mathcal{E}$ or $\dist(E_0,\,\mathcal{E}\setminus E_0)>0$;
\item[(iii)]
either the interior $\textup{int}(E_0)$ of $E_0$ is empty or the boundary of each open component of
$\textup{int}(E_0)$ contains an NCS point.
\end{itemize}
Let $V$ be an open set containing $E_0$ such that $\dist(V,\,\mathcal{E}\setminus E_0)>0$
if $E_0\neq \mathcal{E}$. Then for any asymptotically extremal sequence of monic polynomials
$\{P_n\}_{n\in\mathcal{N}}$ for $E$, 
\begin{equation}\label{eq:main}
\nu_{P_n}|_{V}\,{\stackrel{\star}{\longrightarrow}}\, \mu_E|_{E_0},\quad n\to\infty,\quad n\in\mathcal{N},
\end{equation}
where $\mu|_{K}$ denotes the restriction of a measure $\mu$ to the set $K$.
\end{theorem}

We remark that, for the case of a Jordan region, the hypothesis of Theorem 3 of \cite{GaPo02} implies the existence of an NCS point.
We also note that the assumption $\dist(E_0,\,\mathcal{E}\setminus E_0)>0$ implies that any (open)
component of $\mathrm{int}(E_0)$ is simply connected.\

As a consequence of Theorem \ref{th:main} and \cite[Theorem III.4.1]{ST} we have the following.

\begin{corollary}
With the hypotheses of Theorem~2.1, if $G$ denotes a component of $\mathrm{int}(E_0)$,
then for any asymptotically extremal sequence $\{P_n\}_{n\in\mathcal{N}}$
of monic polynomials for $E$, there exists a point $\zeta$ in $G$ such that
\begin{equation}
\limsup_{n\to\infty,\, n\in\mathcal{N}}|P_n(\zeta)|^{1/n}=\mathrm{cap}(E).
\end{equation}
\end{corollary}

\begin{corollary}\label{cor2.2}
Let $E$ consist of the union of a finite number of  closed Jordan regions $E:=\cup_{j=1}^N \overline{G_j}$, where
$\overline{G_i}\cap\overline{G_j}=\emptyset$, $i\neq j$,
and assume that for each $k=1,\ldots,m$ the boundary of $G_k$ contains an NCS point. Then for any asymptotically extremal sequence of monic polynomials $\{P_n\}_{n\in\mathcal{N}}$ for $E$, 
\begin{equation}
\nu_{P_n}|_{\mathcal{V}}\,{\stackrel{*}{\longrightarrow}}\, \mu_E|_{\mathcal{V}},\quad n\to\infty,\quad n\in\mathcal{N},
\end{equation}
where $\mathcal{V}$ is an open set containing $\bigcup_{k=1}^m\overline{G}_k$, such that if $m<N$ the distance of $\overline{\mathcal{V}}$ from
$\bigcup_{j=m+1}^N \overline{G}_j$ is positive.
\end{corollary}

\begin{figure}
\begin{center}
\includegraphics*[scale=0.8]{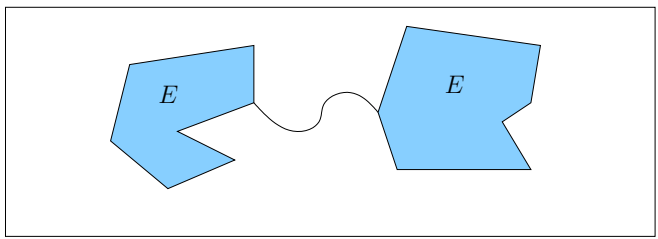}
\end{center}
\caption{Form (iii)}
\label{fig:Form-iii}
\end{figure}

\begin{figure}[h]
\begin{center}\includegraphics*[scale=0.8]{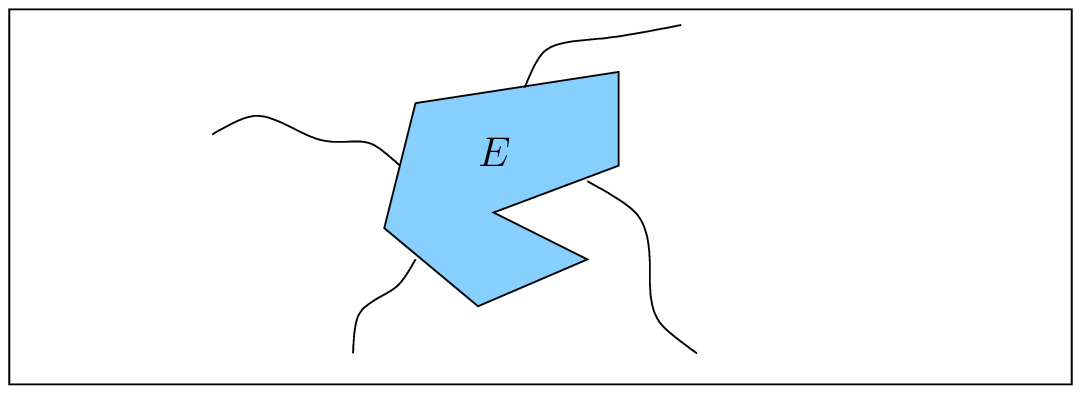}
\end{center}
\caption{Form (iv)}
\label{fig:Form-iv}
\end{figure}

We now give some examples that follow from Theorem~2.1 and Corollary 2.2. If $E$ has one of the following forms, then for any asymptotically extremal sequence $\{P_n\}_{n\in\mathbb{N}}$
of monic polynomials on $E$, we have
\begin{equation}\label{eq:nuPntomuE}
\nu_{P_n}\,{\stackrel{*}{\longrightarrow}}\, \mu_E,\quad n\to\infty,\quad n\in\mathbb{N}.
\end{equation}
\begin{itemize}
\item[(i)]
$E$ is a non-convex polygon or a finite union of mutually exterior non-convex polygons.
\item[(ii)]
$E$ is the union of two mutually exterior non-convex polygons, except for a single common boundary point.
\item[(iii)]
$E$ is the union of two mutually exterior non-convex polygons joined by a Jordan arc in their exterior, such that the complement of $E$ is connected and does not separate the plane; see Figure~\ref{fig:Form-iii}.
\item[(iv)]
$E$ is a a non-convex polygon $\Pi$ together with a finite number of closed Jordan arcs lying exterior to $\Pi$ except for the initial  point on the boundary of $\Pi$,  \textmd{and} such that the complement of $E$ does not separate the plane;
see Figure~\ref{fig:Form-iv}.
\item[(v)]
$E$ is any of the preceding forms with the polygons replaced by closed bounded Jordan regions, each one having an NCS point.
\item[(vi)]
$E$ is any of the preceding forms union with a compact set $K$ in the complement of $E$  such that $K$ has empty interior and $E\cup K$ does not separate the plane.
\end{itemize} 

We remark that if $E$ is a convex region, so that the hypotheses of Theorem \ref{th:main} are not fulfilled, then the zero
behavior of an asymptotically extremal sequence of monic polynomials $P_n$ can be quite different. For example, if
$E$ is the closed unit disk centered at the origin for which $\mu_E$ is the uniform measure on the
 circumference $|z|=1$,  the polynomials $P_n(z)=z^n$ form an extremal sequence for which
$\nu_{P_n}=\delta_0$, the unit point mass at zero. A less trivial example is illustrated
in Figure 3, where the zeros of orthonormal polynomials with respect to area measure on a
circular sector $E$ with opening angle $\pi/2$ are plotted for degrees $n=50, 100,$ and $150.$  These so-called Bergman polynomials $B_n(z)$ form an asymptotically extremal sequence of polynomials  for the sector, yet their normalized zero counting measures converge weak* to a measure
$\nu$ that is supported on the union of three curves lying in the interior of $E$ except for their three
endpoints, the vertices of the sector; see \cite{M-DSS}.\
\begin{figure}[h]
\begin{center}\includegraphics*[scale=0.3]{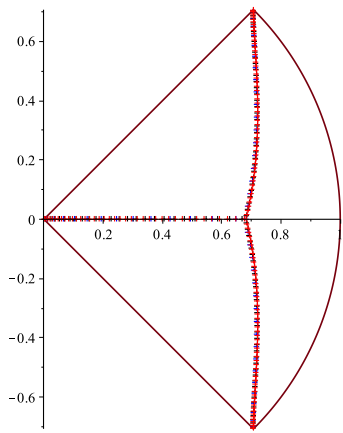}
\end{center}
\caption{Zeros of the Bergman polynomials $B_n$, $n=50,100,150$, for the circular sector with opening angle $\pi/2$.}
\label{cirsec04}
\end{figure}
 On the other hand, for any compact set $E$ of positive capacity, whether convex or not, if $\{q_n(z)\}_{n\in\mathbb{N}}$ denotes a sequence of Fekete polynomials for $E$, then this  sequence
is asymptotically extremal on $E$, all their zeros stay on the outer boundary $\partial\Omega$, and $\nu_{q_n}\,{\stackrel{*}{\longrightarrow}}\, \mu_E,\,$ as $n\to\infty,\,\, n\in\mathbb{N}$; see e.g., \cite [p.~176]{ST}.

In every case, according to Theorem 1.1, a limit measure of a sequence of asymptotically extremal monic polynomials must
have a balayage to the outer boundary of $E$ that equals the equilibrium measure $\mu_E.$  The question then of what types
of point sets can support a measure with such a balayage is a relevant inverse problem. In this connection, there is a conjecture of the first author on the
existence of \textit{electrostatic skeletons} for every convex polygon (more generally, for any convex region with boundary consisting of line segments
or circular arcs). By an electrostatic skeleton on $E$ we mean a positive measure with closed support in $E$, such that its
logarithmic potential  matches the equilibrium potential in $\Omega$ and its support has empty interior and does not separate
the plane. For example, a square region has a skeleton whose support is the union of its diagonals; the circular sector in
Figure 3 has  a  skeleton supported on the illustrated curve joining the three vertices. See the discussion in \cite[p. 55]{LunTo} and in \cite{ErLuRa}.


\section{Applications to special polynomial sequences}
We begin with some results for Bergman polynomials $\{B_n\}_{n\in\mathbb{N}}$ that are orthogonal with respect to the area
$dA$ measure over a bounded Jordan domain $G$; i.e.,
\begin{equation}
\int_G B_m(z)\overline{B_n(z)}dA(z)=0,\quad m\neq n.
\end{equation}

The following theorem was established in \cite{LSS}.
\begin{theorem}\label{LSS-th}
Let $G$ be a bounded Jordan domain whose boundary $\partial G$ is singular; i.e., a conformal map $\varphi$ of $G$ onto the unit disk $\mathbb{D}$ cannot be analytically continued to some open set containing $\overline{G}$. Then, there is a subsequence $\mathcal{N}$ of $\mathbb{N}$ such that
\begin{equation}\label{eq:nuBntomu3.2}
\nu_{B_n}\,{\stackrel{*}{\longrightarrow}}\, \mu_{\partial G},\quad n\to\infty,\quad n\in\mathcal{N}.
\end{equation}
\end{theorem}
It is not difficult to see that this property of $G$ is independent of the choice of the conformal map $\varphi$.

As a consequence of Corollary \ref{cor2.2}, we obtain a result that holds for $\mathcal{N}=\mathbb{N}$.
\begin{corollary}\label{LSS-cor}
If the Jordan domain $G$ has a point on its boundary that satisfies the NCS condition, then (\ref{eq:nuBntomu3.2}) holds for $\mathcal{N}=\mathbb{N}$.
\end{corollary}

\begin{figure}
\begin{center}
\includegraphics*[scale=0.3]{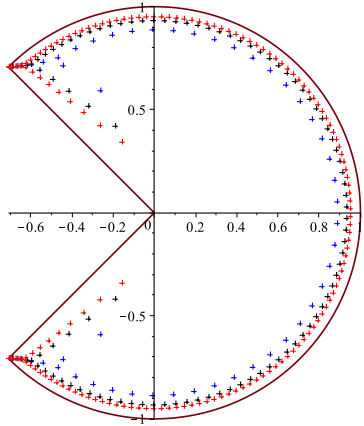}
\end{center}
\caption{Zeros of the Bergman polynomials $B_n$, $n=50,100,150$, for the circular sector with opening angle $3\pi/2$.}
\label{cirsec02}
\end{figure}

In Figure~\ref{cirsec02}, we depict zeros of the Bergman polynomials for the circular sector
$G:=\{z:z=e^{i\theta},\, -3\pi/4<\theta<3\pi/4\}$.
The computations of the Bergman polynomials for this sector as well as for the sector
in Figure 3 were carried out in Maple 16 with 300 significant figures, using the
Arnoldi Gram-Schmidt algorithm; see \cite[Section 7.4]{St-CA13} for a discussion regarding the
stability of the algorithm.

For Figure~\ref{cirsec02} the origin is an NCS point and, therefore, Corollary~\ref{LSS-cor} implies that
the only limit distribution of the zeros is the equilibrium measure, a fact reflected in Figures~\ref{cirsec02}. It
is interesting to note that this figure also depicts the facts that the density of the equilibrium measure  is
zero at corner points with opening angle greater than $\pi$, and it is infinite at corners with opening angle
less that $\pi$.

%

In \cite[p.~1427]{GPSS} the following question has been raised:
If $G$ is the union of two mutually exterior Jordan domains $G_1$ and $G_2$ whose boundaries are singular, does there
exist a common sequence of integers $n$ for which $\nu_{B_n}|_{\mathcal{V}_j}$ converges to
$\mu_{\partial G}|_{\partial G_j}$, where ${\mathcal{V}_j}$ is an open set containing $G_j$, $j=1,2$.
Thanks to Theorem~2.1, the answer is affirmative for the full sequence $\mathbb{N}$
if both $G_1$ and $G_2$ have the NCS property.

We conclude by applying the results of the main theorem to the case of Faber polynomials.
For this, we assume that $\Omega$ is simply connected and let
$\Phi$ denote the conformal map  $\Omega\to\Delta:=\{w:|w|>1\}$, normalized so that near infinity
\begin{equation}\label{eq:Phi}
\Phi(z)=\gamma z+\gamma_0+\frac{\gamma_1}{z}+\frac{\gamma_2}{z^2}+\cdots,\quad \gamma>0.
\end{equation}
The Faber polynomials $\{F_n(z)\}_{n=0}^\infty$ of $E$ are defined as the polynomial part of the expansion of $\Phi^n(z)$
near infinity.

The following theorem was established in \cite{KS}.
\begin{theorem}\label{KS}
Suppose that $\mathrm{int}(E)$ is connected and $\partial E$ is a piecewise analytic curve that has a singularity other
than an outward cusp.
Then, there is a subsequence $\mathcal{N}$ of $\mathbb{N}$ such that
\begin{equation}\label{eq:nuBntomuGcalN}
\nu_{F_n}\,{\stackrel{*}{\longrightarrow}}\, \mu_{E},\quad n\to\infty,\quad n\in\mathcal{N}.
\end{equation}
\end{theorem}

Using Theorem~\ref{th:main} we can refine the above as follows; see also the question raised in Remark 6.1(c) in
\cite{KS}.
\begin{corollary}\label{KS}
If $\mathrm{int}(E)$  has a point on its boundary that satisfies the NCS condition, then (\ref{eq:nuBntomuGcalN})
holds for $\mathcal{N}=\mathbb{N}$.
\end{corollary}

\section{Proof of Theorem~\ref{th:main}}\label{sec:proofs}
Let $\mu$ be any weak$^*$ limit measure of the sequence $\{\nu_{P_n}\}_{n\in\mathcal{N}}$ and recall from Theorem~\ref{th:ST4.7}
that $\supp(\mu)\subset\mathbb{C}\setminus\Omega$ and
\begin{equation}\label{eq:Um=UmE}
U^\mu(z)=U^{\mu_E}(z),\quad  z\in\Omega,
\end{equation}
where
\begin{equation*}
U^\nu(z):=\int\log\frac{1}{|z-t|}d\nu(t)
\end{equation*}
denotes the logarithmic potential on a measure $\nu$.

We consider first the case when $E_0=\mathcal{E}$. It suffices to show that
\begin{equation}\label{eq:supp-mu}
\supp(\mu)\subset\partial\Omega,
\end{equation}
because, this in view of  (\ref{eq:Um=UmE}) and Carleson's unicity theorem (\cite[Theorem II.4.13]{ST})
will imply the relation $\mu=\mu_E$, which yields (\ref{eq:main}) with ${V}=\mathbb{C}$. Clearly, (\ref{eq:supp-mu}) is satisfied automatically in the case $\textup{int}(E_0)=\emptyset$, so we turn our attention now to the case $\textup{int}(E_0)\neq\emptyset$ and assume to the contrary
that $\supp(\mu)$ is not contained in $\partial\Omega$. Then
there exists a small closed disk $K$ belonging
to some open component of $\textup{int}(E_0)$, such that $\mu(K)>0$.
We call this particular component $G$, and note that it is simply connected.
Since $\partial G$ is  regular with respect to the interior
and exterior Dirichlet problem,
\begin{equation*}
\lim_{z\to t\in\partial G \atop{z\in G}}g_G(z,\zeta)=0\mbox{ and }
\lim_{z\to t\in\partial G \atop
{z\in\overline{\mathbb{C}}\setminus\overline{G}}}g_{\overline{\mathbb{C}}\setminus\overline{G}}(z,\infty)=0,
\end{equation*}
where $g_{\overline{\mathbb{C}}\setminus\overline{G}}(z,\infty)$ denotes the Green function of
$\overline{\mathbb{C}}\setminus\overline{G}$ with pole at infinity.

Following Gardiner and Pommeremke (see \cite[Section 5]{GaPo02}), we set $\mu_0:=\mu|_K$ and consider the function
\begin{equation}\label{eq:S-def}
S(z):=\left\{
\begin{array}{cl}
\int g_G(z,\zeta)d\mu_0(\zeta),  &z\in G,\\
U^{\mu_E}(z)-\log\frac{1}{\cpct(E)},   &z\in \overline{\mathbb{C}}\setminus{G}.
\end{array}
\right.
\end{equation}
From the properties of Green functions and equilibrium potentials it follows that $S(z)$ is harmonic in $G\setminus K$ and
in $\Omega\setminus\{\infty\}$, positive in $G$, negative in $\overline{\mathbb{C}}\setminus\overline{G}$ and vanishes
\emph{quasi-everywhere} (that is, apart from a set of capacity zero) on $\partial G$.

Let now $\widehat{\mu}_0$ denote the balayage of $\mu_0$ out of $K$ onto $\partial G$. Then, the relation
$\widehat{\mu}_0\le \mu^b$ follows
from the discussion regarding balayage onto arbitrary compact sets  in \cite[pp. 222--228]{La72}.
Since, by Theorem \ref{th:ST4.7},  $\mu^b=\mu_E$, the difference $\mu_E-\widehat{\mu}_0$ is
a positive measure, a fact leading to the following useful representation:
\begin{equation}\label{eq:S-def-2}
S(z)=U^{\mu_0}(z)+U^{\mu_E-\widehat{\mu}_0}(z)-\log\frac{1}{\cpct(E)}, \quad z\in\mathbb{C},
\end{equation}
which shows that $S(z)$ is superharmonic in $\mathbb{C}$.

By assumption, the boundary of $G$ contains an NCS point $z_0$.  Without loss of generality, we make the following simplifications  regarding the two conditions in
Definition 2.1: By performing a translation and scaling we take $z_0$  to be the origin and by rotation we take
$\zeta_0=i\gamma$, for some $\gamma\in (0,1)$. Finally, in view of the Remark following Definition \ref{Def1}, we take
$\overline{D}:=\{z:|z-i\gamma|\le\gamma\}$ and choose $\gamma$ so that $\overline{D}$ is a subset of the open unit disk $\mathbb{D}$
and $\overline{D}\cap K=\emptyset$.

The contradiction we seek will be a consequence of the following two claims:
\begin{equation*}\label{claim-a}
\textbf{Claim (a)} \quad \int_{\pi}^{2\pi}\frac{S(re^{i\theta})}{r}d\theta \ge \tau,\ as\ r\to 0+,
\end{equation*}
where $\tau$ is a negative constant and
\begin{equation*}\label{claim-b}
\textbf{Claim (b)} \quad \int_{0}^{\pi}\frac{S(re^{i\theta})}{r}d\theta \to +\infty,\ as\ r\to 0+.
\end{equation*}
These claims follow as in \cite{GaPo02}, utilizing in the justification of Claim (b) the essential condition that the
origin is an NCS point so that
%


\begin{equation}\label{eq:limGy}
\lim_{y\to 0+}\frac{g_G(iy,i\gamma)}{y}=+\infty.
\end{equation}\

Note that for small positive $y$ the definition of $S(z)$ gives
\begin{equation}\label{eq:Siy}
\frac{S(iy)}{y}=\frac{1}{y}\int g_G(iy,\zeta)d\mu_0(\zeta).
\end{equation}
Using (\ref{eq:Har}) we have $g_G(iy,\zeta)\ge C g_G(iy,i\gamma)$ for any $\zeta\in K$ and small $y$,  which in view of (\ref{eq:limGy})
leads to the limit
\begin{equation}\label{eq:limSy}
\lim_{y\to 0+}\frac{S(iy)}{y}=+\infty.
\end{equation}

Using Claims (a) and (b) and the fact that $S(0)=0$ (since the origin is a regular point of $\Omega$), it is easy to arrive at a relation that contradicts the mean value inequality for superharmonic
functions (see also \cite[p. 425]{GaPo02}):
\begin{equation*}
\begin{alignedat}{1}
\frac{1}{r}\left(\frac{1}{2\pi}\int_0^{2\pi}S(re^{i\theta})d\theta-S(0)\right)&=
\frac{1}{2\pi r}\int_0^{2\pi}S(re^{i\theta})d\theta\\
&\ge\frac{1}{2\pi r}\int_0^{\pi}S(re^{i\theta})d\theta+\tau\\
&\to\infty, \quad r\to +0.
\end{alignedat}
\end{equation*}
This establishes the theorem for the case $E_0=\mathcal{E}$.

To conclude the proof, we observe that  when $\dist(E_0,\,\mathcal{E}\setminus E_0)>0$ our arguments above show
that $\mu$ cannot have any point of its support in the interior of any open component of ${E_0}\cap V$; hence it is supported on the
outer boundary of $E_0$ inside $V$. Therefore, by following the proof
of Theorem II.4.13 in \cite{ST}, we see that the logarithmic potentials of $\mu$ and $\mu_E$ coincide in $V$ and
the required relation $\mu|_V=\mu_E|_{E_0}$ follows from the unicity theorem for logarithmic potentials.
\qed\\

\noindent\textit{Acknowledgment.} The authors are grateful to the referees for their helpful comments.

\bibliographystyle{amsplain}

\begin{thebibliography}{10}

\bibitem{ArGa}
D.~H. Armitage and S.~J. Gardiner, \emph{Classical {P}otential {T}heory},
  Springer Monographs in Mathematics, Springer-Verlag London, Ltd., London,
  2001.

\bibitem{ErLuRa}
A.~Eremenko, E.~Lundberg, and K.~Ramachandran, \emph{Electrostatic skeletons},
  arXiv (2013).

\bibitem{GaPo02}
S.~J. Gardiner and Ch. Pommerenke, \emph{Balayage properties related to
  rational interpolation}, Constr. Approx. \textbf{18} (2002), no.~3, 417--426.

\bibitem{GPSS}
B.~Gustafsson, M.~Putinar, E.~Saff, and N.~Stylianopoulos, \emph{Bergman
  polynomials on an archipelago: Estimates, zeros and shape reconstruction},
  Advances in Math. \textbf{222} (2009), 1405--1460.

\bibitem{KS}
A.~B.~J. Kuijlaars and E.~B. Saff, \emph{Asymptotic distribution of the zeros
  of {F}aber polynomials}, Math. Proc. Cambridge Philos. Soc. \textbf{118}
  (1995), no.~3, 437--447.

\bibitem{La72}
N.~S. Landkof, \emph{Foundations of {M}odern {P}otential {T}heory},
  Springer-Verlag, New York, 1972, Translated from the Russian by A. P.
  Doohovskoy, Die Grundlehren der mathematischen Wissenschaften, Band 180.

\bibitem{Lehman}
R.~S. Lehman, \emph{Development of the mapping function at an analytic corner},
  Pacific J. Math. \textbf{7} (1957), 1437--1449.

\bibitem{LSS}
A.~L. Levin, E.~B. Saff, and N.~S. Stylianopoulos, \emph{Zero distribution of
  {B}ergman orthogonal polynomials for certain planar domains}, Constr. Approx.
  \textbf{19} (2003), no.~3, 411--435.

\bibitem{LunTo}
E.~Lundberg and V.~Totik, \emph{Lemniscate growth}, Anal. Math. Phys.
  \textbf{3} (2013), no.~1, 45--62.

\bibitem{MhSa91}
H.~N. Mhaskar and E.~B. Saff, \emph{The distribution of zeros of asymptotically
  extremal polynomials}, J. Approx. Theory \textbf{65} (1991), no.~3, 279--300.

\bibitem{M-DSS}
E.~Mina-Diaz, E.~B. Saff, and N.~S. Stylianopoulos, \emph{Zero distributions
  for polynomials orthogonal with weights over certain planar regions}, Comput.
  Methods Funct. Theory \textbf{5} (2005), no.~1, 185--221.

\bibitem{Ra}
T.~Ransford, \emph{Potential {T}heory in the {C}omplex {P}lane}, London
  Mathematical Society Student Texts, vol.~28, Cambridge University Press,
  Cambridge, 1995.

\bibitem{ST}
E.~B. Saff and V.~Totik, \emph{Logarithmic {P}otentials with {E}xternal
  {F}ields}, Springer-Verlag, Berlin, 1997.

\bibitem{StTobo}
H.~Stahl and V.~Totik, \emph{General {O}rthogonal {P}olynomials}, Cambridge
  University Press, Cambridge, 1992.

\bibitem{St-CA13}
N.~Stylianopoulos, \emph{Strong asymptotics for {B}ergman polynomials over
  domains with corners and applications}, Constr. Approx. \textbf{38} (2013),
  no.~1, 59--100.

\end{thebibliography}

\def\cprime{$'$}
\providecommand{\bysame}{\leavevmode\hbox to3em{\hrulefill}\thinspace}
\providecommand{\MR}{\relax\ifhmode\unskip\space\fi MR }
\providecommand{\MRhref}[2]{%
  \href{http://www.ams.org/mathscinet-getitem?mr=#1}{#2}
}
\providecommand{\href}[2]{#2}

\end{document}